\documentclass[a4paper,12pt]{amsart}
\usepackage{amssymb}
\usepackage{ifthen}
\usepackage{graphicx}
\nonstopmode \numberwithin{equation}{section}
\usepackage[usenames]{color}
\newtheorem{thm}{Theorem}
\newtheorem{cor}{Corollary}
\newtheorem{lem}{Lemma}
\newtheorem{prop}{Proposition}
\newtheorem{rem}{Remark}
\newtheorem{rems}[equation]{Remarks}
\theoremstyle{definition}
\newtheorem{defin}{Definition}
\newtheorem{examp}[equation]{Example}
\newtheorem{prob}[equation]{Problem}
\newtheorem{ques}[equation]{Question}
\newtheorem{op}{Problem}

\newtheorem{conj}[equation]{Conjecture}
\newtheorem{deter}[equation]{Determination}
\newtheorem{case}{Case}
\newtheorem{subcase}{Subcase}
\newtheorem{subsubcase}{Subsubcase}
\newtheorem{claim}{Claim}
\newtheorem{subclaim}{Subclaim}

\newcounter {own}
\def\theown {\thesection       .\arabic{own}}

\newenvironment{pf}[1][]{%
 \vskip 3mm
 \noindent
 \ifthenelse{\equal{#1}{}}%
  {{\bf Proof. }}%
  {{\bf #1.} }%
 }%
{\qed\bigskip}

\newcounter{alphabet}
\newcounter{tmp}
\newenvironment{Thm}[1][]{\refstepcounter{alphabet}%
\bigskip%
\noindent%
{\bf Theorem \Alph{alphabet}}%
\ifthenelse{\equal{#1}{}}{}{ (#1)}%
{\bf .} \itshape}{\vskip 8pt}

\makeatletter
\newcommand{\Ref}[1]{\@ifundefined{r@#1}{}{\setcounter{tmp}{\ref{#1}}\Alph{tmp}}}
\makeatother

\newcommand{\IR}{{\mathbb R}}

\def\be{\begin{equation}}
\def\ee{\end{equation}}

\newcommand{\bee}{\begin{enumerate}}
\newcommand{\eee}{\end{enumerate}}

\newcommand{\blem}{\begin{lem}}
\newcommand{\elem}{\end{lem}}
\newcommand{\bthm}{\begin{thm}}
\newcommand{\ethm}{\end{thm}}
\newcommand{\bcor}{\begin{cor}}
\newcommand{\ecor}{\end{cor}}
\newcommand{\beg}{\begin{examp}}
\newcommand{\eeg}{\end{examp}}
\newcommand{\begs}{\begin{examples}}
\newcommand{\eegs}{\end{examples}}
\newcommand{\bdefe}{\begin{defin}}
\newcommand{\edefe}{\end{defin}}
\newcommand{\bprob}{\begin{prob}}
\newcommand{\eprob}{\end{prob}}
\newcommand{\bques}{\begin{ques}}
\newcommand{\eques}{\end{ques}}
\newcommand{\bei}{\begin{itemize}}
\newcommand{\eei}{\end{itemize}}

\newcommand{\bde}{\begin{deter}}
\newcommand{\ede}{\end{deter}}
\newcommand{\bca}{\begin{case}}
\newcommand{\eca}{\end{case}}
\newcommand{\bsca}{\begin{subcase}}
\newcommand{\esca}{\end{subcase}}

\newcommand{\bssca}{\begin{subsubcase}}
\newcommand{\essca}{\end{subsubcase}}

\newcommand{\bcl}{\begin{claim}}
\newcommand{\ecl}{\end{claim}}
\newcommand{\bscl}{\begin{subclaim}}
\newcommand{\escl}{\end{subclaim}}
\newcommand{\bcon}{\begin{conj}}
\newcommand{\econ}{\end{conj}}
\newcommand{\bcons}{\begin{conjs}}
\newcommand{\econs}{\end{conjs}}
\newcommand{\bprop}{\begin{prop}}
\newcommand{\eprop}{\end{prop}}
\newcommand{\br}{\begin{rem}}
\newcommand{\er}{\end{rem}}
\newcommand{\brs}{\begin{rems}}
\newcommand{\ers}{\end{rems}}
\newcommand{\bo}{\begin{obser}}
\newcommand{\eo}{\end{obser}}
\newcommand{\bos}{\begin{obsers}}
\newcommand{\eos}{\end{obsers}}
\newcommand{\bpf}{\begin{pf}}
\newcommand{\epf}{\end{pf}}
\newcommand{\ba}{\begin{array}}
\newcommand{\ea}{\end{array}}
\newcommand{\beq}{\begin{eqnarray}}
\newcommand{\beqq}{\begin{eqnarray*}}
\newcommand{\eeq}{\end{eqnarray}}
\newcommand{\eeqq}{\end{eqnarray*}}

\newcommand{\ds}{\displaystyle}

\newcommand{\bop}{\begin{op}}
\newcommand{\eop}{\end{op}}

\newtheorem{pfofThm1.5}[equation]{}

\newcounter{minutes}
\divide\time by 60
\newcounter{hours}
\multiply\time by 60 \addtocounter{minutes}{-\time}

\begin{document}

\bibliographystyle{amsplain}

\title{On removability properties of $\psi$-uniform domains in Banach spaces}

\def\thefootnote{}
\footnotetext{ \texttt{\tiny File:~\jobname .tex,
          printed: \number\year-\number\month-\number\day,
          \thehours.\ifnum\theminutes<10{0}\fi\theminutes}
} \makeatletter\def\thefootnote{\@arabic\c@footnote}\makeatother

\author{M. Huang}
\address{M. Huang, Department of Mathematics,
Hunan Normal University, Changsha,  Hunan 410081, People's Republic
of China} \email{mzhuang79@163.com}

\author{M. Vuorinen}
\address{M. Vuorinen, Department of Mathematics and Statistics,
University of Turku, 20014 Turku,
Finland
} \email{vuorinen@utu.fi}

\author{X. Wang $^* $
}
\address{X. Wang, Department of Mathematics,
Hunan Normal University, Changsha,  Hunan 410081, People's Republic
of China} \email{xtwang@hunnu.edu.cn}

\date{}
 \keywords{
Uniform domain, $\psi$-uniform domain,
 removability property, quasihyperbolic metric.\\
${}^{\mathbf{*}}$ Corresponding author}

\begin{abstract}
Suppose that $E$ denotes a real Banach space with the
dimension at least $2$. The main aim of this paper is to show that a
domain $D$ in $E$ is a $\psi$-uniform domain if and only if $D\backslash P$ is a
$\psi_1$-uniform domain, and  $D$  is a uniform domain if and only if $D\backslash P$ also is a
uniform domain, whenever $P$ is a closed countable subset of $D$ satisfying
a quasihyperbolic separation condition. This condition requires
 that the quasihyperbolic distance (w.r.t. $D$) between each pair of
distinct points in $P$ has a lower bound greater than or equal to
$\frac{1}{2}$.
\end{abstract}

\thanks{The research was partly supported by NSFs of
China (No. 11071063 and No. 11101138) and by the Academy of Finland,
Project 2600066611.}

\maketitle{} \pagestyle{myheadings} \markboth{}{The removability property of $\psi$-uniform domains and uniform domains in Banach spaces}

\section{Introduction and main results}\label{sec-1}

The quasihyperbolic metric of a domain in a metric space was introduced by F. W. Gehring and
his students B. Palka and B. Osgood in the 1970's \cite{geos,GP} in the setup of the Euclidean space
${\mathbb R}^n, n\ge2.$ Since its first appearance, the quasihyperbolic metric has become
an important tool in geometric function theory and in its
generalizations to metric spaces and to Banach spaces \cite{Vai5}.
For instance, V\"ais\"al\"a's theory of quasiconformal maps
in Banach spaces relies on the quasihyperbolic metric.
Yet, some basic questions of the quasihyperbolic geometry in Banach
spaces and even in Euclidean spaces are open. Such a basic question as the convexity of quasihyperbolic balls has been studied in \cite{k,krt,rt,Vai9}  only recently.

In this paper, we study the classes of uniform domains \cite{masa}
and the wider class of $\psi$-uniform domains  \cite{vu85} in Banach
spaces and the stability of these classes of domains under the
removal of a countable set of points. The motivation for this study
stems from the discussions in \cite{Mw,Vai6}. In these two papers,
the removability questions were studied for the classes of uniform
domains and John domains.  For latest results
in the case of $\psi$-uniform domains in $\mathbb{R}^n$
see \cite{klsv}.  Our main result extends these results to Banach spaces.
We begin with some basic definitions and the statements of our
results. The proofs and necessary supplementary notation and
terminology will be given thereafter.

Throughout the paper, we always assume that $E$ denotes a real
Banach space with the dimension at least $2$. The norm of a
vector $z$ in $E$ is written as $|z|$, and for each pair of points
$z_1$, $z_2$ in $E$, the distance between them is denoted by
$|z_1-z_2|$ and the closed line segment with endpoints $z_1$ and $z_2$
by $[z_1, z_2]$. We always use $\mathbb{B}(x_0,r)$ to denote the
open ball $\{x\in E:\,|x-x_0|<r\}$ centered at $x_0$ with radius
$r>0$. Similarly, for the closed balls and spheres, we employ the
usual notations $\overline{\mathbb{B}}(x_0,r)$ and $\partial
\mathbb{B}(x_0,r)$, respectively. When $x_0=0$, briefly, we denote $\mathbb{B}(x_0,r)=\mathbb{B}(r)$, in particular, $\mathbb{B}=\mathbb{B}(1)$. In $\IR^n$, we specially denote $\mathbb{B}^n(x_0,r)$ the ball centered at
$x_0$ with radius $r$, $\mathbb{B}^n(r)$ the ball centered at
$0$ with radius $r$, and $\mathbb{B}^n$ the ball centered at
$0$ with radius $1$. We adopt some basic terminology
following closely \cite{Martio-80, Vai, Vai6-0, Vai6}.

\bdefe \label{def1} A domain $D$ in $E$  is said to be  {\it
$c$-uniform}  if there exists a constant $c$ with the property that
each pair of points $z_{1},z_{2}$ in $D$ can be joined by a
rectifiable arc $\gamma$ in $ D$ satisfying
\begin{enumerate}
\item\label{eq-1}
$\ds\min_{j=1,2}\ell (\gamma [z_j, z])\leq c\, d_D(z)$ for all $z\in
\gamma$, and

\item\label{eq-2}
$\ell(\gamma)\leq c\,|z_{1}-z_{2}|$,
\end{enumerate}
where $\ell(\gamma)$ denotes the arc length of $\gamma$,
$\gamma[z_{j},z]$ the part of $\gamma$ between $z_{j}$ and $z$, and
$d_D(z)$ the distance from $z$ to the boundary $\partial D$ of $D$ \cite{masa}.
Also we say that $\gamma$ is a {\it double $c$-cone arc}.\edefe

\bdefe \label{def1.2} Let $\psi:[0,\infty)\to [0,\infty)$ be a
homeomorphism. A domain $D$ in $E$ is called  $\psi$-{\it uniform}
if $$k_{D}(z_1,z_2)\leq \psi\Big(\frac{|z_1-z_2|}{\min\{d_{D}(z_1),
d_{D}(z_2)\}}\Big)$$ for all  $z_{1}$, $z_{2}\in D$ \cite{vu85}. \edefe

It is a basic fact that
uniform domains are  special cases of $\psi$-uniform domains.
See Section $2\,.$


Simple examples show, see Example \ref{countnonpsi}, that removing a countable
closed set $E$ from a uniform domain $D$ may yield a domain $D\setminus E$ which
is not $\psi$-uniform for any $\psi$.
This motivates us to study countable
subsets of a domain satisfying the following quasihyperbolic
separation condition.

We say that a countable closed set $P$
in a domain $D\subset E$ satisfies a quasihyperbolic separation condition if
 \begin{equation}
 \label{india-1}  P=\{x_j\in D:   j=1,2,\dots\} \,\,{\rm and}\,\,
k_D(x_i,x_j)\geq \frac{1}{2} \,, \quad  \forall i\not=j\,.
\end{equation}
In what follows, for a given domain $D,$
the symbol $P$ stands for a fixed countable subset of $D$ with this property.

The purpose of this paper is to study the following problem.

\bop \label{Con1}Suppose that $D$ is a $\psi$-uniform (resp. $c$-uniform) domain in $E$
and $P$ is as in \eqref{india-1}. Is it true that
 $G=D\backslash P$ is a $\psi_1$-uniform (resp. $c_1$-uniform) domain, where
 $\psi_1$ (resp. $c_1$) depends only on $\psi$ (resp. $c$)?\eop


We are now in a position to formulate our results.

\begin{thm}\label{th1.1}  A domain $D$ in $E$ is a $\psi$-uniform domain if and
only if $G=D\backslash P$ is a $\psi_1$-uniform domain, where
$\psi$ and $\psi_1$ are $(2^{12}, 1, 3, 2^7)$-equivalent.\end{thm}

Here for two homeomorphisms $\psi$ and $\psi_1:\; [0,\infty)\to [0,\infty)$, we say that they are
 $(a_1, a_2, a_3, a_4)$-{\it equivalent} if there are positive constants $a_1$, $a_2$, $a_3$ and $a_4$
 such that $\psi_1(t)=a_1\psi(a_2t)$ and $\psi(t)=a_3\psi_1(a_4t)$
 for $t\geq0$.
As a corollary of Theorem \ref{th1.1}, we have

\bcor\label{th1.2} A domain $D$ in $E$ is a $c$-uniform domain if and
only if $G=D\backslash P$ is a $c_1$-uniform domain, where the
constants $c$ and $c_1$ depend only on each other.\ecor

\br
It easily follows from the discussions in Section \ref{sec-4} that Theorem \ref{th1.1} and Corollary \ref{th1.2}
still hold when we replace the constant $\frac{1}{2}$ in $P$ in \eqref{india-1} by any positive constant $\kappa$.
\er

The proofs of Theorem \ref{th1.1} and Corollary \ref{th1.2} will be given in Section \ref{sec-4}. We will
 prove several lemmas in Section \ref{sec-3}, which will be used later on, and in Section \ref{sec-2}, some preliminaries will
be introduced.

\section{Preliminary results}\label{sec-2}


\subsection{Quasihyperbolic distance and neargeodesics}
\ \ \

The {\it quasihyperbolic length} of a rectifiable arc or a path
$\alpha$ in the norm metric in $D$ is the number:

$$\ell_{k_D}(\alpha)=\int_{\alpha}\frac{|dz|}{d_{D}(z)}.
$$

Gehring and Palka \cite{GP} introduced the quasihyperbolic metric of
a domain in $\IR^n$ and it has been recently used by many authors
 in the study of quasiconformal mappings \cite{Ahl,HK,Vai9,vu85,vu88} and related questions \cite{HIMPS}.
Many of the
  basic properties of this metric may be found in \cite{geos,Vai6-0, Vai6, Vai3,vu88}.

For each pair of points $z_1$, $z_2$ in $D$, the {\it distance ratio
metric} $j_D(z_1,z_2)$ between $z_1$ and $z_2$ is defined by
$$j_D(z_1,z_2)=\log\Big(1+\frac{|z_1-z_2|}{\min\{d_D(z_1),d_D(z_2)\}}\Big).$$

 The {\it
quasihyperbolic distance} $k_D(z_1,z_2)$ between $z_1$ and $z_2$ is
defined in the usual way:
$$k_D(z_1,z_2)=\inf\{\ell_{k_D}(\alpha)\},
$$
where the infimum is taken over all rectifiable arcs $\alpha$
joining $z_1$ to $z_2$ in $D$. For all $z_1$, $z_2$ in $D$, we have
\cite{Vai3}

\beq\label{eq(0000)} k_{D}(z_1, z_2)\geq
\inf_{\alpha}\left\{\log\Big(1+\frac{\ell(\alpha)}{\min\{d_{D}(z_1),
d_{D}(z_2)\}}\Big)\right\}\geq \Big|\log
\frac{d_{D}(z_2)}{d_{D}(z_1)}\Big|,\eeq where the infimum is taken
over all rectifiable curves $\alpha$ in $D$ connecting $z_1$ and
$z_2$. Since $\ell(\alpha)\ge |z_1-z_2|$ in \eqref{eq(0000)}, for all $z_1$, $z_2$ in $D$, we have
\be\label{vvm-1}k_D(z_1,z_2)\geq j_D(z_1,z_2).\ee

Next, if $|z_1-z_2|\le d_D(z_1)$, then we have \cite{Vai6-0},
\cite[Lemma 2.11]{vu81}
\begin{equation} \label{upperbdk}
k_D(z_1,z_2)\le \log\Big( 1+ \frac{
|z_1-z_2|}{d_D(z_1)-|z_1-z_2|}\Big)  \le \frac{
|z_1-z_2|}{d_D(z_1)-|z_1-z_2|}\,,
\end{equation}
where the last inequality follows from the following elementary
inequality
\begin{equation} \label{logINE}
\frac{r}{1-r/2} \le \log \frac{1}{1-r} \le \frac{r}{1-r} \, \quad {\rm  for }\,\, 0\le r<1 \,.
\end{equation}

The following observation easily follows from \eqref{vvm-1} and
Definition \ref{def1.2}.
\bprop\label{prop-1} If $D$ is
$\psi$-uniform, then the homeomorphism $\psi$ satisfies
$$\psi(t)\geq \log(1+t)$$ for $t>0$.\eprop

In \cite{Vai6},  V\"ais\"al\"a characterized uniform domains by the
quasihyperbolic metric.

\begin{Thm}\label{thm0.1} {\rm (\cite{geos, Vai6})}
\;
For a domain $D$, the following are quantitatively equivalent: \bee

\item $D$ is a $c$-uniform domain;
\item $k_D(z_1,z_2)\leq c'\;
 j_D(z_1,z_2)$ for all $z_1,z_2\in D$;
\item $k_D(z_1,z_2)\leq c'_1\;
 j_D(z_1,z_2)+d$ for all $z_1,z_2\in D$.\eee
\end{Thm}

 In the case of domains in $ {\mathbb R}^n \,,$ the equivalence
 of items (1) and (3) in Theorem D is due to Gehring and Osgood \cite{geos} and the
 equivalence of items (2) and (3) due to Vuorinen \cite{vu85}.
By Theorem \Ref{thm0.1}, we see that uniformity implies
$\psi$-uniformity.

\begin{examp} \label{countnonpsi}
For $r>0$ and $j=10,11,...$ let $E(r,j)$ be a finite set of points in
$A=\mathbb{B}^n( \frac{21\,r}{20}) \setminus \mathbb{B}^n( r)$ such that for every point $z \in A$
we have $d(z,A)<r/(20j)\,.$ Let $G(r,j)= {\mathbb R}^n \setminus E(r,j)$ and
$a_j = \frac{11\,r}{10}e_1\,,$ $b_j = \frac{9\,r}{10}e_1\,.$ It is readily seen that
$$k_{G(r,j)}(a_j,b_j)\ge \frac{r/20}{r/(20\, j)} = j\,.$$
Let $E= \cup_{j=10}^{\infty} E(2^{-j},j) \cup \{ 0\} \,.$ Clearly, $\mathbb{B}^n \setminus E$
is not $\psi$-uniform for any $\psi$ although $\mathbb{B}^n$ is uniform.
\end{examp}


Recall that an arc $\alpha$ from $z_1$ to $z_2$ is a {\it
quasihyperbolic geodesic} if $\ell_{k_D}(\alpha)=k_D(z_1,z_2)$. Each
subarc of a quasihyperbolic geodesic is obviously a quasihyperbolic
geodesic. It is known that a quasihyperbolic geodesic between every
pair of points in $E$ exists if the dimension of $E$ is finite, see
\cite[Lemma 1]{geos}. This is not true in Banach spaces
\cite[Example 2.9]{Vai4}. In order to remedy this shortage,
V\"ais\"al\"a introduced the following concepts \cite{Vai6}.

\bdefe \label{def1.5} Let $D$ be a domain in $E$ and $\nu>1\,.$
An arc $\alpha\subset D$ is a {\it $\nu$-neargeodesic} if
$\ell_{k_D}(\alpha[x,y])\leq \nu\;k_D(x,y)$ for all $x, y\in
\alpha$.\edefe

Obviously, a $\nu$-neargeodesic is a quasihyperbolic geodesic if and
only if $\nu=1$. The smoothness of geodesics has been studied
recently in \cite{RT2}.

In \cite{Vai4}, V\"ais\"al\"a proved the following property
concerning the existence of neargeodesics in Banach spaces.

\begin{Thm}\label{LemA} $($\cite[Theorem 3.3]{Vai4}$)$
Fix  $\nu>1\,.$  Then for
 $\{z_1,\, z_2\}\subset D \subset E$  there is a
$\nu$-neargeodesic in $D$ joining $z_1$ and $z_2$.
\end{Thm}

\subsection{Quasiconvexity}

\bdefe\label{def1.3} We say that an arc $\gamma$ in $D\subset E$ is
$c$-{\it quasiconvex} in the norm metric if it satisfies the
condition
$$\ell(\gamma[z_1,z_2])\leq c\,|z_1-z_2|
$$
for every $z_1$, $z_2$ in $\gamma$.
 \edefe

The following result is  due to Sch\"affer \cite{Sc}.


\begin{Thm}\label{LemC}{\rm (\cite[4.4]{Sc})} Suppose that $S$ in $E$ is a sphere, that
$T$ is a $2$-dimensional linear subspace in $E$ and that the
intersection $S\cap T$ contains at least two points. For every pair
 $\{z_1,\;z_2\}\subset T\cap S$, if $\gamma\subset T\cap S$ is the minor arc or a half circle
 with the endpoints $z_1$ and $z_2$,
 then $\gamma$ is $2$-quasiconvex. \end{Thm}





\section{Several Lemmas}\label{sec-3}

We recall that $D$ denotes a domain in $E$ and $G=D\backslash P$,
where $P\subset D$ is a countable set satisfying the
quasihyperbolic separation condition \eqref{india-1}. The aim of this section
is to prove several lemmas on which the proofs of Theorem \ref{th1.1} and
Corollary \ref{th1.2} will be based.
Our first lemma concerns the number of the points of $P$, which are
contained in the balls  $\mathbb{B}(x,\lambda
d_D(x))$ in $D$ for varying values of $\lambda$.

\begin{lem}\label{lem2.1-0} $(1)$
Let $ \lambda_0 \equiv 1-\exp(-1/4)=0.22\ldots\,.$
If  $
\lambda\in (0, \lambda_0)\,$ and $x\in D\,,$ then $\mathbb{B}(x,\lambda
d_D(x))$ contains at most one point of $P\,.$

$(2)$ If $\lambda\in (0, \lambda_0)\,,$
$x\in G$, and $\frac{d_G(x)}{d_D(x)}<\lambda$, then
$\mathbb{B}(x,\lambda d_D(x))$ contains exactly one point of $P\,.$

$(3)$ If  $\lambda\leq \frac{1}{16}$, $x\in G$, and
$\frac{d_G(x)}{d_D(x)}<\lambda$, then $\mathbb{B}(x,\lambda d_D(x))$
contains exactly one point $x_i$ in $P$ which satisfies
$$d_G(z)=|x_i-z|$$ for all $z\in \overline{\mathbb{B}}(x,\frac{1}{16}d_D(x))$.

\end{lem}

\bpf $(1)$
Let $x\in D$. It follows from (\ref{upperbdk}) that
$$k_D(\mathbb{B}(x,\lambda d_D(x)))=
\sup\{k_D(u,v):u,v\in \mathbb{B}(x,\lambda d_D(x))\}\leq 2\log
\frac{1}{1-\lambda}< \frac{1}{2}.$$
Then \eqref{india-1} implies that $(1)$ holds.

$(2)$ Let $x\in G$ with $\frac{d_G(x)}{d_D(x)}<\lambda$. Then there
exists some $i$ such that $|x-x_i|<\lambda d_D(x)$. Hence $x_i\in
\mathbb{B}(x,\lambda d_D(x))$ and by $(1)$ $\mathbb{B}(x,\lambda
d_D(x))$ cannot contain points in $P\backslash \{x_i\}$.

We now prove $(3)$. Because
$\frac{1}{5} \in (0,\lambda_0),$ we see by $(1)$ that
$\mathbb{B}(x,\frac{1}{5} d_D(x))$ contains at most one point of $P$.
By $(2)$ there exists a unique point
$x_i\in \mathbb{B}(x,\lambda d_D(x)) \cap P\,.$ Then
$$d_G(x)=|x-x_i|< \lambda \, d_D(x),$$
because for $j \neq i$ we have
$$ |x-x_j|\ge \frac{1}{5} d_D(x) -  \lambda \, d_D(x) > \lambda \, d_D(x) $$
by the choice of $\lambda\,.$ In the same way we see that for all
$z \in \mathbb{B}(x,\frac{1}{15} d_D(x)) $ we have for all $j \neq i$
$$
|z-x_j| \ge |x-x_j|- |x-z| >  \frac{1}{5} d_D(x) -  \lambda \, d_D(x)=(\frac{1}{5} -\lambda) d_D(x)
$$
whereas
$$
|z-x_i| \le  |z-x|+ |x-x_i| <  \frac{1}{15} d_D(x) +
\lambda \, d_D(x)=(\frac{1}{15} +\lambda) d_D(x)\,.
$$
Because $\frac{1}{5} - \lambda> \frac{1}{15}  + \lambda$, we see that
for all $z\in
\overline{\mathbb{B}}(x,\frac{1}{16}d_D(x))$,
 $$d_G(z)=|x_i-z|  \,.$$
The proof is complete.\epf

The constant $\lambda_0$ in Lemma \ref{lem2.1-0} obviously depends on the constant $\frac{1}{2}$ in the quasihyperbolic separation condition \eqref{india-1}. At the same time it is easy to see that a result similar to  Lemma \ref{lem2.1-0} also holds if the constant $\frac{1}{2}$ in the definition of the set $P$ \eqref{india-1} is replaced by $\sigma \in (0,\frac{1}{2}) \,.$ Only the constants would
change and e.g. $\lambda_0$ would be replaced by $\lambda_{\sigma}\,,$ with
$$  2\log
\frac{1}{1-\lambda_{\sigma}} = \sigma \,. $$

Given a point $w$ in a domain $U \subset {\mathbb R}^n$ it is clear that
$$   k_U(x,y) \le k_{U \setminus \{ w\}}(x,y)$$
for all $x,y \in {U \setminus \{ w\}\,}.$ For points $x,y \in  U \setminus
{\mathbb B}(w, \theta d_U(w)),\,0< \theta<1\,,\,$ we also have an opposite inequality as
shown in \cite[Lemma 2.53]{vu85}.
The next few lemmas  deal with this situation for Banach spaces.
We start by a comparison theorem for the metrics $k_G$
and $j_G\,.$

\begin{lem}\label{lem2.1-1}
For $0<\mu\leq \frac{1}{32}$ and $w_1$, $w_2\in G$, if $w_2\in
\overline{\mathbb{B}}(w_1, \mu d_D(w_1))$ and $\min\{d_G(w_1),
d_G(w_2)\}\leq \frac{\mu}{2}d_D(w_1)$, then
$$k_G(w_1,w_2)\leq \frac{13}{2}j_G(w_1,w_2).$$
\end{lem}

\bpf Clearly, we have 

 \beq\label{mm-2-3}\max\{d_G(w_1),
d_G(w_2)\}&\leq& \min\{d_G(w_1), d_G(w_2)\}+|w_1-w_2|\\
\nonumber&\leq&\frac{3\mu}{2}d_D(w_1).\eeq

 By Lemma \ref{lem2.1-0} $(3)$, we see that there exists
some point $x_i\in \mathbb{B}(w_1,\frac{3\mu}{2}d_D(w_1))\cap P$
satisfying
$$|w_2-x_i|=d_G(w_2)\;\;\mbox{and}\;|w_1-x_i|=d_G(w_2).$$

Without loss of generality, we may assume
\be\label{mm-2-4}\min\{d_G(w_1), d_G(w_2)\}=d_G(w_2).\ee


We use $w_{1,1}$ to denote the intersection point of the closed
segment $[w_1,x_i]$ with the sphere $\mathbb{S}(x_i, \min\{d_G(w_1),
d_G(w_2)\})$. It is possible that $w_{1,1}=w_1$. Let $T$ denote a
$2$-dimensional linear subspace of $E$ passing thorough the points
$w_1$, $w_2$ and $x_i$, and $\omega_0$ the circle $T\cap
\mathbb{S}(x_i, \min\{d_G(w_1), d_G(w_2)\})$. Then $w_{1,1}$ and
$w_2$ divide the circle $\omega_0$ into two parts $\beta_i$ and
$\beta_{1,i}$. Without loss of generality, we may assume that
$\ell(\beta_i)\leq \ell(\beta_{1,i})$. Then it follows from
(\ref{mm-2-4}) that for each $z\in\beta_i$,  $$|z-w_1|\leq
|w_1-x_i|+|x_i-z|\leq 2\mu d_D(w_1),$$ whence
$$\beta_i\subset \mathbb{\overline{B}}\big(w_1,
2\mu d_D(w_1)\big),$$ and so Lemma \ref{lem2.1-0} yields that for
each $z\in\beta_i$,
$$d_G(z)=d_G(w_2),$$
which, together with (\ref{upperbdk}) and Theorem \Ref{LemC}, shows
that
\begin{eqnarray*}k_G(w_1,w_2)&\leq& k_G(w_1,w_{1,1})+\ell_{k_G}(\beta_i)
\\
\nonumber&\leq&
\log\Big(1+\frac{|w_1-w_{1,1}|}{d_G(w_2)}\Big)+\frac{2|w_2-w_{1,1}|}{d_G(w_2)}\\
\nonumber&\leq& \log\Big(1+\frac{|w_1-w_{1,1}|}{d_G(w_2)}\Big)+ \frac{4}{\log 3}\log \Big(1+\frac{|w_2-w_{1,1}|}{d_G(w_2)}\Big)\\
\nonumber&<&  (1+\frac{6}{\log 3})\log\Big(1+
\frac{|w_1-w_2|}{d_G(w_2)}\Big)\\
\nonumber&\leq& \frac{13}{2}j_G(w_1,w_2),\end{eqnarray*} since
$\frac{|w_2-w_{1,1}|}{d_G(w_2)}\leq 4$ and $|w_2-w_{1,1}|\leq 2
|w_1-w_2|$. Hence the proof follows.\epf

The next two results are related to $k_G$ and $k_D$.

\begin{lem}\label{lem2.1-1-1}
For $w_1\in G$, suppose $d_G(w_1)= \frac{1}{128}d_D(w_1)$. If
$w_2\in \mathbb{S}(w_1, \frac{1}{32}d_D(w_1))$, then
$$k_G(w_1,w_2)\leq 2^9k_D(w_1,w_2).$$
\end{lem}
\bpf It follows from Lemma \ref{lem2.1-0} that there is a unique
element $x_i$ in the intersection
$P\cap\overline{\mathbb{B}}(w_1,\frac{1}{64}d_D(w_1))$ such that
$$d_G(w_2)=|w_2-x_i|\geq \frac{3}{128}d_D(w_1)=3d_G(w_1).$$
Then Lemma \ref{lem2.1-1}, \eqref{eq(0000)} and the Bernoulli inequality \cite[(3.6)]{vu88}
imply
\begin{eqnarray*}k_G(w_1,w_2)&\leq& \frac{13}{2}\log\Big(1+
\frac{|w_1-w_2|}{d_G(w_1)}\Big)\\
\nonumber&=&\frac{13}{2}\log\Big(1+
\frac{128|w_1-w_2|}{d_D(w_1)}\Big)\\ \nonumber&<& 2^9\log\Big(1+
\frac{|w_1-w_2|}{d_D(w_1)}\Big)\\
\nonumber&\leq& 2^9k_D(w_1,w_2),\end{eqnarray*} which shows that the
lemma is true.\epf

\begin{lem}\label{lem2.1-2}
Let $w_1$, $w_2\in G$ and let $\gamma$ denote a $2$-neargeodesic joining
$w_1$ and $w_2$ in $D$. If $d_G(z)\geq \frac{1}{128}d_D(z)$ for each
$z\in\gamma$, then
$$k_G(w_1,w_2)\leq 2^8k_D(w_1,w_2).$$
\end{lem}

\bpf Obviously, we get
 \begin{eqnarray*}k_G(w_1,w_2)&\leq&
\ell_{k_G}(\gamma[w_1,w_2]) =
\int_{\gamma[w_1,w_2]}\frac{|dx|}{d_G(x)}
\leq
128\int_{\gamma[w_1,w_2]}\frac{|dx|}{d_D(x)}\\
\nonumber&=&128\ell_{k_D}(\gamma[w_1,w_2])\leq 2^8k_D(w_1,w_2).\end{eqnarray*} Hence the proof is
complete.\epf

Our last lemma in this section is as follows.

\begin{lem}\label{lem2.1-3}
Let $w_1$, $w_2\in D$. If $|w_1-w_2|\geq \frac{1}{c}d_D(w_1)$
$(c\geq 2)$, then
$$|w_1-w_2|\geq \frac{1}{c+1}d_D(w_2).$$
\end{lem}

\bpf Suppose on the contrary that $$|w_1-w_2|< \frac{1}{c+1}d_D(w_2).$$
Then
$$d_D(w_1)\geq d_D(w_2)-|w_1-w_2|>\frac{c}{c+1}d_D(w_2),$$ which shows that
$$|w_1-w_2|\geq \frac{1}{c}d_D(w_1)>\frac{1}{c+1}d_D(w_2).$$ This is the desired contradiction.\epf

\section{The proofs of Theorem \ref{th1.1} and Corollary \ref{th1.2}}\label{sec-4}

\subsection{The proof of Theorem \ref{th1.1}}
We first prove the sufficiency. Suppose that $G$ is a
$\psi_1$-uniform domain. Then we shall prove that for  $z_1$,
$z_2\in D$,
\beq\label{india-2} k_D(z_1,z_2)\leq
\psi\Big(\frac{|z_1-z_2|}{\min\{d_D(z_1),d_D(z_2)\}}\Big),\eeq where
$\psi(t)=3\psi_1(2^7t)$ for $t>0$, which implies that $D$ is a
$\psi$-uniform domain. Without loss of generality, we assume that
$$\min\{d_G(z_1), d_G(z_2)\}=d_G(z_1).$$

We divide the proof into two cases.

\subsubsection{ We first suppose that $\;|z_1-z_2|\leq \frac{1}{2}d_D(z_1)$} \label{xvv-1}\ \

Then it follows from  (\ref{upperbdk}) and Proposition \ref{prop-1} that
\beq\label{sun-1}
k_D(z_1,z_2) &\leq &
\log\Big(1+\frac{|z_1-z_2|}{d_D(z_1)-|z_1-z_2|}\Big)\leq
2\log\Big(1+\frac{|z_1-z_2|}{d_D(z_1)}\Big)\\ \nonumber   &\leq & 2\psi_1\Big(\frac{|z_1-z_2|}{\min\{d_D(z_1),d_D(z_2)\}}\Big).\eeq

\subsubsection{ We then suppose that $\;|z_1-z_2|> \frac{1}{2}d_D(z_1)$}\label{xvv-1}\ \

Hence by Lemma \ref{lem2.1-3}, we have \be\label{mzxt-1}|z_1-z_2|\geq
\frac{1}{3}d_D(z_2).\ee

In the following, we separate the rest discussions to two parts.

\bca $d_G(z_1)\geq \frac{1}{64}d_D(z_1).$\eca
Under this assumption, we have
\beq\label{sun-2} k_D(z_1,z_2)&\leq & k_G(z_1,z_2)\leq
\psi_1\Big(\frac{|z_1-z_2|}{d_G(z_1)}\Big)\leq
\psi_1\Big(\frac{64|z_1-z_2|}{d_D(z_1)}\Big)\\ \nonumber &\leq & \psi_1\Big(\frac{64|z_1-z_2|}{\min\{d_D(z_1),d_D(z_2)\}}\Big).\eeq

\bca $d_G(z_1)< \frac{1}{64}d_D(z_1).$\eca
For a proof under this assumption, we let $u_1\in
\mathbb{S}(z_1,\frac{1}{32}d_D(z_1))$. Then by \eqref{upperbdk}
 \be\label{ccc-0}k_D(z_1,u_1)\leq
\frac{32|z_1-u_1|}{31d_D(z_1)}=\frac{1}{31}<
\log\Big(1+\frac{|z_1-z_2|}{d_D(z_1)}\Big)\,.\ee
Moreover, we get \be\label{ccc-3}|u_1-z_2|\leq
|z_1-z_2|+|z_1-u_1|\leq \frac{17}{16}|z_1-z_2|,\ee and it follows
from Lemma \ref{lem2.1-0} and the assumption ``$d_G(z_1)<
\frac{1}{64}d_D(z_1)$" that there exists only one element $x_i$ in
$\overline{\mathbb{B}}(z_1,\frac{1}{64}d_D(z_1))\cap P$ such that
$d_G(u_1)=|u_1-x_i|$ and so \be\label{ccc-2} d_D(u_1)>
d_G(u_1)=|u_1-x_i|\geq |u_1-z_1|-|z_1-x_i|\geq
\frac{1}{64}d_D(z_1).\ee

In this part, we again distinguish two possibilities.

\bsca $d_G(z_2)\geq \frac{1}{120}d_D(z_2).$\esca
Then by (\ref{ccc-0}), (\ref{ccc-3}) and (\ref{ccc-2}), we have
\beq\label{sun-3} k_D(z_1,z_2)&\leq&k_D(z_1,u_1)+k_D(u_1,z_2)\\
\nonumber&\leq&
\log\Big(1+\frac{|z_1-z_2|}{d_D(z_1)}\Big)+k_G(u_1,z_2)
\\ \nonumber&\leq& j_D(z_1,z_2)+\psi_1\Big(\frac{|u_1-z_2|}{\min\{d_G(u_1),d_G(z_2)\}}\Big)
\\ \nonumber&\leq& j_D(z_1,z_2)+\psi_1\Big(\frac{120|u_1-z_2|}{\min\{d_D(z_1),d_D(z_2)\}}\Big)
\\ \nonumber&\leq& j_D(z_1,z_2)+\psi_1\Big(\frac{2^7|z_1-z_2|}{\min\{d_D(z_1),d_D(z_2)\}}\Big)
\\ \nonumber&<&
2\psi_1\Big(\frac{2^7|z_1-z_2|}{\min\{d_D(z_1),d_D(z_2)\}}\Big).\eeq

\bsca $d_G(z_2)< \frac{1}{120}d_D(z_2).$\esca
Under this hypothesis, to get a homeomorphism $\psi$ from $\psi_1$, we take $u_2\in
\mathbb{S}(z_2,\frac{1}{32}d_D(z_2))$. It follows from Lemma
\ref{lem2.1-0} and (\ref{mzxt-1}) that there is only one element
$x_i$ in $\overline{\mathbb{B}}(z_2,\frac{1}{32}d_D(z_2))\cap P$
such that $d_G(u_2)=|u_2-x_i|$ and so \be\label{ccc-4} d_D(u_2)\geq
d_G(u_2)=|u_2-x_i|\geq |u_2-z_2|-|z_2-x_i|\geq
\frac{11}{480}d_D(z_2)\ee and by \eqref{mzxt-1}, we have
\be\label{ccc-5}|u_1-u_2|\leq |z_1-z_2|+|z_1-u_1|+|u_2-z_2|\leq
\frac{37}{32}|z_1-z_2|.\ee

It follows from \eqref{mzxt-1}, \eqref{upperbdk} and
(\ref{ccc-0}) that
\be\label{ccc-6}k_D(z_2,u_2)\leq
\frac{32|z_2-u_2|}{31d_D(z_2)}<
\log\Big(1+\frac{|z_1-z_2|}{d_D(z_2)}\Big).\ee
Then we infer from
\eqref{ccc-0}, (\ref{ccc-2}), (\ref{ccc-4}), (\ref{ccc-5}),
(\ref{ccc-6})  and Proposition \ref{prop-1} that \beq\label{sun-4}
k_D(z_1,z_2)&\leq&k_D(z_1,u_1)+k_G(u_1,u_2)+k_D(u_2,z_2)\\
\nonumber&\leq& 2j_D(z_1,z_2)+k_G(u_1,u_2)
\\
\nonumber&\leq& 2j_D(z_1,z_2)
+\psi_1\Big(\frac{|u_1-u_2|}{\min\{d_G(u_1),d_G(u_2)\}}\Big)\\
\nonumber&\leq& 2j_D(z_1,z_2)
+\psi_1\Big(64\frac{|u_1-u_2|}{\min\{d_D(z_1),d_D(z_2)\}}\Big)\\
\nonumber&\leq& 2j_D(z_1,z_2)
+\psi_1\Big(74\frac{|z_1-z_2|}{\min\{d_D(z_1),d_D(z_2)\}}\Big)\\
\nonumber&<&
3\psi_1\Big(74\frac{|z_1-z_2|}{\min\{d_D(z_1),d_D(z_2)\}}\Big).\eeq
Since $\psi_1$ is increasing, by taking $\psi(t)=3\psi_1(2^7t)$ for $t\geq 0$, we easily see
 from the inequalities
\eqref{sun-1}, \eqref{sun-2}, \eqref{sun-3} and \eqref{sun-4} that \eqref{india-2} holds.
\medskip

Next we prove the necessity. Suppose that $D$ is a
$\psi$-uniform domain. Then we shall prove that for $z_1$, $z_2\in
G$,
\beq\label{india-3}
k_G(z_1,z_2)\leq 2^{12}\psi\Big(\frac{|z_1-z_2|}{\min\{d_G(z_1),d_G(z_2)\}}\Big),\eeq which implies that
$G$ is a $\psi_1$-uniform domain with $\psi_1=2^{12}\psi$.

Without loss of generality, we may assume that
$\min\{d_G(z_1),d_G(z_2)\}=d_G(z_1)$. In the following, we consider
the two cases where $d_{G}(z_1)\leq \frac{1}{64}d_D(z_1)$ and
$d_{G}(z_1)> \frac{1}{64}d_D(z_1)$, respectively.

\subsubsection{We first suppose that $\;d_{G}(z_1)\leq
\frac{1}{64}d_D(z_1)$}\label{Hvw-2-1}\ \
\medskip

Let $\gamma$ be a $2$-neargeodesic joining $z_1$ and $z_2$ in $D$. The existence of $\gamma$
follows from Theorem \Ref{LemA}. We separate the discussions in this case to two parts.

\bca\label{mvw-1} $|z_1-z_2|\leq \frac{1}{32}d_D(z_1)$.\eca

Then by Lemma \ref{lem2.1-1} and Proposition \ref{prop-1}, we have
\bcl\label{cl-1} $k_G(z_1,z_2)\leq \frac{13}{2}\log\Big(1+
\frac{|z_1-z_2|}{d_G(z_1)}\Big)\leq
\frac{13}{2}\psi\Big(\frac{|z_1-z_2|}{d_G(z_1)}\Big).$\ecl

 \bca\label{mvw-2}
$|z_1-z_2|> \frac{1}{32}d_D(z_1)$.\eca Under this condition, Lemma \ref{lem2.1-3}
implies \beq\label{sun-5}|z_1-z_2|\geq \frac{1}{33}d_D(z_2).\eeq
 Obviously, there exists
some point $v_1\in \gamma\cap \mathbb{S}(z_1, \frac{1}{32}d_D(z_1))$
such that
$$\gamma[z_2,v_1]\subset D\backslash\mathbb{B}(z_1,
\frac{1}{32}d_D(z_1)).$$ By Lemma \ref{lem2.1-0}, there exists some
point $x_{i,1}\in \overline{\mathbb{B}}(z_1,\frac{1}{64}d_D(z_1))\cap
P$ such that \be\label{hv-3-1-0}d_G(v_1)=|v_1-x_{i,1}|\geq
\frac{1}{64}d_D(z_1)\geq \frac{1}{66}d_D(v_1),\ee since
$d_D(v_1)\leq d_D(z_1)+|z_1-v_1|\leq \frac{33}{32}d_D(z_1)$.

It follows from Lemma \ref{lem2.1-1} and (\ref{hv-3-1-0}) that
\be\label{hv-3-1-1}k_G(z_1,v_1)\leq\frac{13}{2}\log\Big(1+
\frac{|z_1-v_1|}{d_G(z_1)}\Big)\leq\frac{13}{2}\log\Big(1+
\frac{|z_1-z_2|}{d_G(z_1)}\Big).\ee

In this part, we again distinguish two possibilities.

\bsca \label{mxm-1}$d_G(z)\geq \frac{1}{128}d_D(z)$ for each
$z\in\gamma[v_1,z_2]$.\esca
By Lemma \ref{lem2.1-2}, we know
\begin{eqnarray*} k_G(v_1,z_2)&\leq&
2^8k_D(v_1,z_2)\leq2^8\ell_{k_D}(\gamma[v_1,z_2])\leq2^8\ell_{k_D}(\gamma[z_1,z_2])\\
\nonumber &\leq& 2^9k_D(z_1,z_2)\leq
2^9\psi\Big(\frac{|z_1-z_2|}{d_G(z_1)}\Big),\end{eqnarray*} since
$D$ is $\psi$-uniform, and so the following inequality easily
follows from (\ref{hv-3-1-1}). \bcl\label{mxm-x-1} $k_G(z_1,z_2)\leq
k_G(z_1,v_1)+k_G(v_1,z_2)\leq
519\psi\Big(\frac{|z_1-z_2|}{d_G(z_1)}\Big).$\ecl

\bsca \label{mxm-2} There exists some point
$z\in\gamma[v_1,z_2]$ such that $d_G(z)<
\frac{1}{128}d_D(z)$.\esca

Under this assumption, it follows from (\ref{hv-3-1-0}) that there exists some point $y_1$
which is the first point in $\gamma$ along the direction from $v_1$
to $z_2$ such that $$d_G(y_1)=\frac{1}{128}d_D(y_1).$$ Then Lemma
\ref{lem2.1-2} shows \beq\label{mzz-1}k_G(v_1,y_1)\leq
2^8k_D(v_1,y_1) \leq2^8\ell_{k_D}(\gamma[z_1,z_2])\\
\nonumber \leq 2^9k_D(z_1,z_2)\leq
2^9\psi\Big(\frac{|z_1-z_2|}{d_D(z_1)}\Big).\eeq

To get a homeomorphism $\psi_1$ from $\psi$ in this possibility, we consider two cases.

\bssca\label{sat-4} $|z_2-y_1|\leq\frac{1}{32}d_D(y_1).$\essca
Then we see from Lemma \ref{lem2.1-1} that
\be\label{mmm-1}k_G(y_1,z_2)\leq
\frac{13}{2}\log\Big(1+\frac{|y_1-z_2|}{\min\{d_G(y_1),d_G(z_2)\}}\Big).\ee
Thus we have the following claim.

\bcl\label{mz-1-1} $k_G(z_1,z_2)\leq
564\psi\Big(\frac{|z_1-z_2|}{d_G(z_1)}\Big).$\ecl

We now prove this claim. Since $$d_D(z_2)\geq
d_D(y_1)-|z_2-y_1|\geq \frac{31}{32}d_D(y_1)\geq 31|z_2-y_1|,$$ we
infer from \eqref{sun-5} that \beq\label{mz-03} |z_1-z_2|\geq
\frac{31}{33}|y_1-z_2|,\eeq and by \eqref{mmm-1}, we get
$$k_G(y_1,z_2)\leq
\frac{13}{2}\log\Big(1+\frac{33|z_1-z_2|}{31\min\{d_G(y_1),d_G(z_2)\}}\Big).$$
Let us leave the proof of Claim \ref{mz-1-1} for a moment and prove the following inequality
  \beq\label{mz-1} k_G(y_1,z_2)\leq
35\log\Big(1+\frac{|z_1-z_2|}{d_G(z_2)}\Big).\eeq

To prove this estimate, obviously, we only need to consider the case
$d_G(z_2)\geq d_G(y_1)$. Since $$d_G(z_2)\leq d_G(y_1)+|z_2-y_1|\;\;\mbox{and}\;\;
|z_2-y_1|\leq \frac{1}{32}d_D(y_1)= 4d_G(y_1),$$
we see from \eqref{mmm-1} and the Bernoulli inequality \cite[(3.6)]{vu88} that
$$k_G(y_1,z_2)\leq
\frac{13}{2}\log\Big(1+\frac{|y_1-z_2|}{d_G(y_1)}\Big)\leq
\frac{13}{2}\log\Big(1+\frac{5|y_1-z_2|}{d_G(z_2)}\Big)<
35\log\Big(1+\frac{|z_1-z_2|}{d_G(z_2)}\Big).$$ Hence \eqref{mz-1}
is true.\medskip

Now, we come back to the proof of Claim \ref{mz-1-1}. It follows from
(\ref{hv-3-1-1}), (\ref{mzz-1}) and \eqref{mz-1} that
\begin{eqnarray*}k_G(z_1,z_2)&\leq&
k_G(z_1,v_1)+k_G(v_1,y_1)+k_G(y_1,z_2)\\ \nonumber&\leq&
\frac{13}{2}\log\Big(1+
\frac{|z_1-z_2|}{d_G(z_1)}\Big)+2^9\psi\Big(\frac{|z_1-z_2|}{d_G(z_1)}\Big)+35\log\Big(1+\frac{|z_1-z_2|}{d_G(z_2)}\Big)
\\ \nonumber&\leq&42\log\Big(1+\frac{|z_1-z_2|}{d_G(z_1)}\Big)+2^9\psi\Big(\frac{|z_1-z_2|}{d_G(z_1)}\Big)
\\
\nonumber&\leq&564\psi\Big(\frac{|z_1-z_2|}{d_G(z_1)}\Big),\end{eqnarray*}
which shows that Claim \ref{mz-1-1} is true.

\bssca\label{s-s-2} $|z_2-y_1|>\frac{1}{32}d_D(y_1)$.\essca Obviously,
there exists some point $v_2\in \gamma\cap \mathbb{S}(y_1,
\frac{1}{32}d_D(y_1))$ such that $$\gamma[z_2,v_2]\subset
D\backslash\mathbb{B}(y_1, \frac{1}{32}d_D(y_1)).$$

By Lemma \ref{lem2.1-0}, we see that there exists some point
$x_{i,2}\in P\cap\mathbb{\overline{B}}(y_1,\frac{1}{128}d_D(y_1))$
such that \beq\label{hv-3-1-01}d_G(v_2)&=&|v_2-x_{i,2}|\geq
|v_2-y_1|-|y_1-x_{i,2}|\\ \nonumber&\geq&
\frac{3}{128}d_D(y_1)\geq\frac{1}{44}d_D(v_2),\eeq since
$d_D(v_2)\leq d_D(y_1)+|y_1-v_2|\leq \frac{33}{32}d_D(y_1)$.

Moreover, by Lemma \ref{lem2.1-1-1}, we have
\be\label{mmm-3-1}k_G(y_1,v_2)\leq 2^9k_D(y_1,v_2)\leq
2^{10}k_D(z_1,z_2).\ee

If $d_G(z)\geq \frac{1}{128}d_D(z)$ for each $z\in\gamma[v_2,z_2]$,
then Lemma \ref{lem2.1-2} implies \beq\label{sat-1} k_G(v_2,z_2)\leq
2^8k_D(v_2,z_2)\leq 2^9k_D(z_1,z_2),\eeq and thus we infer from
\eqref{hv-3-1-1},
 \eqref{mzz-1}, \eqref{mmm-3-1}, \eqref{sat-1} that
$$k_G(z_1,z_2)\leq \frac{13}{2}\log\Big(1+\frac{|z_1-z_2|}{d_G(z_1)}\Big)+ 3\cdot2^{10}k_D(z_1,z_2).$$ Hence we have the
following estimate.
 \bcl
$k_G(z_1,z_2)\leq 2^{12}\psi\Big(\frac{|z_1-z_2|}{d_G(z_1)}\Big).$\ecl

For the remaining case, that is, there exists some point
$z\in\gamma[v_2,z_2]$ such that $d_G(z)< \frac{1}{128}d_D(z)$,
similar arguments as in Subcase \ref{mxm-1} show that there exists
some point $y_2\in \gamma[v_2,z_2]$ satisfying
$$d_G(y_2)=\frac{1}{128}d_D(y_2)$$ and \beq \label{sat-2} k_G(v_2,
y_2)\leq 2^8k_D(v_2, y_2).\eeq Now, if $|z_2-y_2|\leq
\frac{1}{32}d_D(y_2)$, then  the similar reasoning  as in the proof
of \eqref{mz-1} shows that \beq\label{sat-3} k_G(y_2,z_2)\leq
35\log\Big(1+\frac{|z_1-z_2|}{d_G(z_2)}\Big).\eeq Then it follows
from \eqref{hv-3-1-1}, \eqref{mzz-1}, \eqref{mmm-3-1}, \eqref{sat-2}
and \eqref{sat-3} that \begin{eqnarray*}k_G(z_1, z_2)&\leq& k_G(z_1,
v_1)+ k_G(v_1, y_1)+k_G(y_1, v_2) +
 k_G(v_2, y_2)+ k_G(y_2, z_2) \\ \nonumber&\leq&  2^{10}k_D(z_1,z_2)+42\log\Big(1+\frac{|z_1-z_2|}{d_G(z_1)}\Big).\end{eqnarray*}
 Hence we
reach the following estimate.

 \bcl $k_G(z_1,
z_2)\leq 2^{12}\psi\Big(\frac{|z_1-z_2|}{d_G(z_1)}\Big).$\ecl

We assume now that $|z_2-y_2|> \frac{1}{32}d_D(y_2)$. Then there
exists some point $v_3\in \gamma[y_2,z_2]\cap
\mathbb{S}(y_2,\frac{1}{32}d_D(y_2))$ such that
$$\gamma[z_2,v_3]\subset
D\backslash\mathbb{B}(y_2, \frac{1}{32}d_D(y_2)),$$ and the similar
reasoning as in the proof of (\ref{mmm-3-1}) shows that
$$k_G(y_2,v_3)\leq 2^9k_D(y_2,v_3)\leq
2^{10}k_D(z_1,z_2).$$

By repeating the procedure as above, we will reach a finite sequence
of points in $\gamma$: \begin{enumerate}

\item\label{tue-1}   $\{z_1, v_1,y_1,v_2,y_2,\cdots, v_t, z_2\}$ such that
$d_G(z)\geq \frac{1}{128}d_D(z)$ for each $z\in \gamma[v_t, z_2]$;
or

\item\label{tue-2}
 $\{z_1, v_1,y_1,v_2,y_2,\cdots, v_t, y_t, z_2\}$ such that
 $|z_2-y_t|\leq \frac{1}{32}d_D(y_t)$.\end{enumerate}

It follows from \eqref{eq(0000)} that
$$k_D(y_i,v_{i+1})\geq \log\Big(1+\frac{|y_i-v_{i+1}|}{d_D(y_i)}\Big)\geq \log\frac{33}{32}$$ for each
$i\in\{1,\cdots,t\}$, and we see that $$t\leq
\frac{k_D(z_1,z_2)}{\log\frac{33}{32}}.$$

For the former case, i.e., when the statement \eqref{tue-1} as above
holds, we have shown that
\begin{enumerate}
\item $k_G(z_1,v_1)\leq\frac{13}{2}\log\Big(1+
\frac{|z_1-v_1|}{d_G(z_1)}\Big)\leq\frac{13}{2}\log\Big(1+
\frac{|z_1-z_2|}{d_G(z_1)}\Big)$;
\item $k_G(v_i,
y_i)\leq 2^8k_D(v_i, y_i)$, where $i\in \{1, \cdots, t-1\}$;
\item $k_G(y_i,v_{i+1})\leq 2^9k_D(y_i,v_{i+1})$, where $i\in \{1, \cdots, t\}$; and
\item $k_G(v_t, z_2)\leq 2^8k_D(v_t,z_2)\leq 2^9k_D(z_1,z_2)$. \end{enumerate}
Hence we obtain
\begin{eqnarray*}k_G(z_1,z_2)&\leq&  k_G(z_1,v_1)+\sum_{i=1}^{t-1} k_G(v_i,y_i)+\sum_{i=2}^{t} k_G(y_{i-1}, v_i)+k_G(v_t,z_2)
\\ \nonumber&\leq& \frac{13}{2}\log\Big(1+ \frac{|z_1-z_2|}{d_G(z_1)}\Big)+
2^8\sum_{i=1}^{t-1} k_D(v_i,y_i)+ 2^9\sum_{i=2}^{t}
k_D(y_{i-1},v_i)\\
\nonumber&&+35\log\Big(1+\frac{|z_1-z_2|}{d_G(z_2)}\Big)
\\ \nonumber&\leq& \frac{13}{2}\log\Big(1+ \frac{|z_1-z_2|}{d_G(z_1)}\Big)+
2^{10}k_D(z_1,z_2)+35\log\Big(1+\frac{|z_1-z_2|}{d_G(z_2)}\Big),
\end{eqnarray*} which
shows
\bcl\label{sun-6} $k_G(z_1,z_2)\leq
2^{11}\psi\Big(\frac{|z_1-z_2|}{d_G(z_1)}\Big)$.\ecl

For the latter case, i.e., the statement \eqref{tue-2} as above
holds, we also have shown that
\begin{enumerate}
\item $k_G(z_1,v_1)\leq\frac{13}{2}\log\Big(1+
\frac{|z_1-v_1|}{d_G(z_1)}\Big)\leq\frac{13}{2}\log\Big(1+
\frac{|z_1-z_2|}{d_G(z_1)}\Big)$;
\item $k_G(v_i,
y_i)\leq 2^8k_D(v_i, y_i)$, where $i\in \{1, \cdots, t\}$;
\item $k_G(y_i,v_{i+1})\leq 2^9k_D(y_i,v_{i+1})$, where $i\in \{1, \cdots, t-1\}$, and
\item $k_G(y_t,z_2)\leq
35\log\Big(1+\frac{|z_1-z_2|}{d_G(z_2)}\Big)$. \end{enumerate} Hence
we get
\begin{eqnarray*}k_G(z_1,z_2)&\leq&  k_G(z_1,v_1)+\sum_{i=1}^{t} k_G(v_i,y_i)+\sum_{i=2}^{t} k_G(y_{i-1}, v_i)+k_G(y_t,z_2)
\\ \nonumber&\leq& \frac{13}{2}\log\Big(1+ \frac{|z_1-z_2|}{d_G(z_1)}\Big)+
2^8\sum_{i=1}^{t} k_D(v_i,y_i)+ 2^9\sum_{i=2}^{t}
k_D(y_{i-1},v_i)\\
\nonumber&&+35\log\Big(1+\frac{|z_1-z_2|}{d_G(z_2)}\Big)
\\ \nonumber&\leq& \frac{13}{2}\log\Big(1+ \frac{|z_1-z_2|}{d_G(z_1)}\Big)+
2^{10}k_D(z_1,z_2)+35\log\Big(1+\frac{|z_1-z_2|}{d_G(z_2)}\Big),
\end{eqnarray*} which
implies \bcl\label{sun-7} $k_G(z_1,z_2)\leq
2^{11}\psi\Big(\frac{|z_1-z_2|}{d_G(z_1)}\Big)$.\ecl By taking
$\psi_1(t)=2^{12}\psi(t)$ for $t\geq 0$, we see that in the case $d_G(z_1)\leq
\frac{1}{64}d_D(z_1)$, \eqref{india-3} follows from Claims \ref{cl-1} $\sim$
\ref{sun-7}.

\subsubsection{We then suppose that $\;d_G(z_1)>
\frac{1}{64}d_D(z_1)$}\label{Hvw-2-2}\ \
\medskip

Also, we divide the discussions into two cases.

\bca $|z_1-z_2|\leq \frac{1}{2}d_G(z_1)$.\eca
By \eqref{upperbdk} we have
$$k_G(z_1,z_2)\leq
\frac{2|z_1-z_2|}{d_G(z_1)}\,.$$
Hence by Proposition \ref{prop-1}, we have
\bcl\label{sun-10}
$k_G(z_1,z_2)\leq
3\log\Big(1+\frac{|z_1-z_2|}{d_G(z_1)}\Big)\leq 3\psi\Big(\frac{|z_1-z_2|}{d_G(z_1)}\Big).$
\ecl

\bca $|z_1-z_2|> \frac{1}{2}d_G(z_1)$\eca
Since $|z_1-z_2|> \frac{1}{128}d_D(z_1)$, we know from
Lemma \ref{lem2.1-3} that
\be\label{ccc-1} |z_1-z_2|>\frac{1}{129}d_D(z_2).\ee

Let $\beta$ be a $2$-neargeodesic joining $z_1$ and $z_2$ in $D$.
We divide the rest discussions into two subcases.
 \bsca $d_G(z)\geq
\frac{1}{64}d_D(z)$ for each $z\in\beta$.\esca In this case, the
following inequality easily follows from Lemma \ref{lem2.1-2}.
\bcl\label{sun-11} $k_G(z_1,z_2)\leq2^8k_D(z_1,z_2)\leq
2^8\psi\Big(\frac{|z_1-z_2|}{d_G(z_1)}\Big).$\ecl

\bsca There exists some point $z\in\beta$ such that $d_G(z)<
\frac{1}{64}d_D(z)$.\esca Obviously, it follows from the assumption
``$d_G(z_1)> \frac{1}{64}d_D(z_1)$" that there exists point $p_1$
which is the first point in $\beta$ along the direction from $z_1$
to $z_2$ such that
$$d_G(p_1)=\frac{1}{64}d_D(p_1).$$
Then Lemma \ref{lem2.1-2} shows \beq\label{tue-3} k_G(z_1, p_1) \leq
2^8k_D(z_1, p_1).\eeq

To get a homeomorphism $\psi_1$ from $\psi$, we consider the case where $|z_2-p_1|\leq
\frac{1}{32}d_D(p_1)$ and the case where
$|z_2-p_1|>\frac{1}{32}d_D(p_1)$, respectively.

\bssca $|z_2-p_1|\leq \frac{1}{32}d_D(p_1)$. \essca

Under this condition, it follows from (\ref{ccc-1}) that \be\label{cctv-1}|z_2-p_1|\leq
\frac{1}{31}d_D(z_2)\leq \frac{129}{31}|z_1-z_2|\ee since
$d_D(z_2)\geq d_D(p_1)-|z_2-p_1|\geq \frac{31}{32}d_D(p_1)$.

By Lemma \ref{lem2.1-0}, we have
$$ d_G(z_2)\leq \frac{1}{16}d_D(p_1)\leq 4d_G(p_1).$$
Then we know from Lemma \ref{lem2.1-1} and (\ref{cctv-1}) that
$$k_G(p_1, z_2)\leq \frac{13}{2}\log\Big(1+
\frac{|z_2-p_1|}{\min\{d_G(z_2),
d_G(p_1)\}}\Big)\leq28\log\Big(1+\frac{|z_2-z_1|}{d_G(z_2)}\Big),$$
which, together with \eqref{tue-3}, implies \bcl\label{sun-12}
$k_G(z_1, z_2)\leq 28\psi\Big(\frac{|z_1-z_2|}{d_G(z_1)}\Big)$. \ecl
\bssca $|z_2-p_1|> \frac{1}{32}d_D(p_1)$. \essca Obviously, there
exists some point $q_1\in \beta[p_1,z_2]$ such that
$$\beta[q_1,z_2]\subset
D\backslash\mathbb{B}(p_1, \frac{1}{32}d_D(p_1)).$$
By Lemma \ref{lem2.1-0}, we see that there exists some point
$x_{i,3}\in P\cap\mathbb{\overline{B}}(p_1,\frac{1}{128}d_D(p_1))$
such that \beq\label{hv-3-1-01}d_G(q_1)&=&|q_1-x_{i,3}|\geq
|q_1-p_1|-|p_1-x_{i,3}|\\ \nonumber&\geq&
\frac{3}{128}d_D(p_1)\geq\frac{1}{44}d_D(q_1),\eeq since
$d_D(q_1)\leq d_D(p_1)+|p_1-q_1|\leq \frac{33}{32}d_D(p_1)$.

Then the similar reasoning as in Subsubcase \ref{s-s-2} in Subsection
\ref{Hvw-2-1} implies that we will get a finite sequence of points in $\beta$:

\begin{enumerate}
\item   $\{z_1, p_1, q_1, \cdots, p_s, z_2\}$ such that
$d_G(z)\geq \frac{1}{128}d_D(z)$ for each $z\in \gamma[p_s, z_2]$;
or

\item
 $\{z_1, p_1,q_1,\cdots, p_s, q_s, z_2\}$ such that
 $|z_2-q_s|\leq \frac{1}{32}d_D(q_s)$.
\end{enumerate}
It follows from the similar arguments as in Claims \ref{sun-6} and
\ref{sun-7} in Subsubsection \ref{s-s-2} that \bcl\label{sun-13}
$k_G(z_1, z_2)\leq 2^{11}\psi\Big(\frac{|z_1-z_2|}{d_G(z_1)}\Big)$.
\ecl
Now we are in a
position to conclude that the truth of \eqref{india-3} for the case $d_G(z_1)>
\frac{1}{64}d_D(z_1)$ follows from Claims \ref{sun-10} $\sim$ \ref{sun-13}.
Hence the proof of
Theorem \ref{th1.1} is complete.\qed

\subsection{The proof of Corollary \ref{th1.2}}
First, we prove the sufficiency. Suppose $G=D\backslash P$ is a
$c_1$-uniform domain. Then Theorem \Ref{thm0.1} implies that there
exists a constant $c_1'$ depending only on $c_1$ such that for all
$x$ and $y$ in $G$,
$$k_G(x,y)\leq c'_1\;
 j_G(x,y).$$
By Theorem \ref{th1.1}, we see that
 for all $z_1$ and $z_2$ in $D$,
$$k_D(z_1, z_2)\leq 3c'_1\log\Big(1+
\frac{128|z_1-z_2|}{\min\{d_D(z_1),d_D(z_2)\}}\Big)\leq 384c'_1
j_D(z_1, z_2),$$ which, together with Theorem \Ref{thm0.1}, shows
that $D$ is a $c$-uniform domain, where $c$ depends only on $c_1$.

Next, we prove the necessity. Suppose that $D$ is a $c$-uniform domain.
Then Theorem \Ref{thm0.1} implies that there exists a constant $c'$
depending only on $c$ such that for all $x$ and $y$ in $D$,
$$k_D(x,y)\leq c'\;
 j_D(x,y).$$
By Theorem \ref{th1.1}, we see that
 for all $z_1$ and $z_2$ in $G=D\backslash P$,
$$k_G(z_1, z_2)\leq 2^{12}c'j_G(z_1, z_2),$$ which, together
with Theorem \Ref{thm0.1},  shows that $G$ is a $c_1$-uniform
domain, where $c_1$ depends only on $c$. \qed

\bigskip

\bigskip

{\bf Acknowledgement.} This research was carried out when the first author was an academic visitor
 at the University of Turku and supported by the Academy of Finland grant of the second author with the
Project number 2600066611.  The authors are indebted to the referee for his/her
valuable suggestions.

\end{document}